\documentclass[12pt]{article}
\usepackage{amsthm}
\usepackage{amsfonts}
\usepackage{geometry}
\usepackage{graphicx}     
\usepackage{svg}

\usepackage{hyperref}     
\newtheorem{theorem}{Theorem}[section]
\newtheorem{lemma}[theorem]{Lemma}

\newtheorem{conjecture}[theorem]{Conjecture}
\newtheorem{definition}[theorem]{Definition}

\newtheorem{remark}[theorem]{Remark}

\begin{document}

\title{Knotted Surfaces, Homological Norm and Extendable Subgroup\footnote{This article is published on 24 Oct 2025 by  {\it Top. and its App.}. See \href{https://doi.org/10.1016/j.topol.2025.109644}{https://doi.org/10.1016/j.topol.2025.109644}}}
\author{Qiling Liu \\Peking University\\ 2201110018@pku.edu.cn}
\date{\today}
\maketitle

\begin{abstract}
\noindent We prove that for an arbitrary $g$, there is a surface $K$ of genus $g$ embedded in $S^4$, which has finitely many extendable self-homeomorphisms' action on $H_1(K,\mathbb{Z})$, by defining a norm on $H_1(K,\mathbb{Z})$ and proving its additivity.
\end{abstract}

\section{Introduction}
Assume that there is a $p$-manifold smoothly embedded in $S^{p+2}$ via $i\colon M \to S^{p+2}$. We would like to consider the problem as to whether and how many self-homeomorphisms of $M$ can extend to the whole space $S^{p+2}$. For a knotted $T^p$ in $S^{p+2}$ or an orientable closed surface $F_g$ of genus $g$ in $S^4$, this has been studied by many authors \cite{oo}\cite{spin}. For an orientation preserving self-homeomorphism $h\colon M\to M$, we call $h$ is extendable over $i$ if there is an orientation preserving self-homeomorphism of $S^{p+2}$, $\tilde{h}\colon S^{p+2}\to S^{p+2}$, such that its limit on $M$ is $h$, namely $i \circ h=\tilde{h} \circ i$. All extendable self-homeomorphisms of $i$ form a group, which we denote by $E(i)$.

Therefore, the problem is that, when given an embedding $i$ what can we say about $E(i)$? Some results have completely calculated $E(i)$ in some simple cases. For the trivial embedding $i$ from $F_g$ to $S^4$, article \cite{triv} gives a principle to judge if an element belongs to $E(i)$. Whether there is an $F_g$ embedded in $S^4$ having no non-trivial extendable self-homeomorphism is still open. In this article, we will give an embedding whose $E(i)$ has finite image in $\mathrm{Aut}(F_g,\mathbb{Z})=\mathrm{Sp}(2g,\mathbb{Z})$, namely the image of mapping class group $\mathrm{MCG}(F_g)$'s action on $H_1(F_g,\mathbb{Z})$(with some symplectic basis), by norm method.

A simple case is when $M=T^p$ and $i\colon T^p\to S^{p+2}$ is standardly unknotted in $S^{p+2}$. Article \cite{attr} analyses kinds of extendable diffeomorphisms and proves the index of $E(i)$ in $\mathrm{Aut}(T^p)$, $[\mathrm{Aut}(T^p):E(i)]$ is at most $2^p-1$, where $\mathrm{Aut}(T^p)$ is the linear automorphism group of $H_1(T^p,\mathbb{Z})$, namely $\mathrm{SL}(p,\mathbb{Z})$; and by another result as shown in \cite{spin} we can obtain $[\mathrm{Aut}(T^p):E(i)]$ is exactly $2^p-1$. Article \cite{spin} also proves that $[\mathrm{MCG}(T^p):E(i)]$ is finite in differential and PL category. 

A more complicated case is $T^p$ knotted in $\mathbb{R}^{p+2}$. Article \cite{spin} uses the spin method and gives estimate of $[\mathrm{MCG}(T^p):E(i)]$. Article \cite{triv} determines the group $E(i)$ for trivial $F_g$ embedded in $S^4$, and calculates $[\mathrm{MCG}(F_g):E(i)]=2^{2g-2}+2^{g-1}$ as actually the lower bound \cite{spin} obtains. Also, in articles \cite{hiro} and \cite{surg}, the authors determine the group $E(i)$ for some cases of $T^2$ knotted in $S^4$.

Another result is in article \cite{oo}, where the authors define a new norm similar to the Thurston norm, but on $H_1(T^2)$. The authors use this method to prove that there is some embedding $i\colon T^2\to S^{4}$ such that $E(i)$ is finite. This shows knotted case is quite different from the unknotted case. A good question is whether this is true for dimension $p>2$. In this article, we care about the knotted case of high genus surface $F_g$ and corresponding version of the theorem.

In another sight, for a surface $F_g$, one can fix some element $f\in \mathrm{MCG}(F_g)$ and discuss whether there exists an embedding such that $f$ is extendable. In \cite{spin}, the authors prove for any $g\geq 1$, there exists $f\in \mathrm{MCG}_{top}(F_g)$, which is not homeomorphically extendable over any smooth embedding $i\colon F_g\to \mathbb{R}^4$. 

Also see \cite{peri} for more details of finite order homeomorphisms.

In this article, for a surface $M=F_g$, denote the mapping class group's action on $H_1(M,\mathbb{Z})$(with some symplectic basis) by $\mathrm{Aut}(M,\mathbb{Z})=\mathrm{Sp}(2g,\mathbb{Z})$. We prove that 

\begin{theorem}\label{main}
For any $g\geq 1$, there is a surface $M=F_g$ of genus g embedded in $S^4$ such that the image of its extendable self-homeomorphisms in $\mathrm{Aut}(M,\mathbb{Z})=\mathrm{Sp}(2g,\mathbb{Z})$(with some symplectic basis) can only be diagonal matrixes with diagonal elements $\pm1$. 
\end{theorem}

Compared with the previous results, this theorem extends them to the case of  ``genus bigger than 1''. In order to prove the theorem, we continue to use the norm method as shown in \cite{oo}. 

This article is organized as follows. We will introduce the norm's definition in Section 2, and prove the additivity of connected sum in Section 3. In Section 4 we apply the norm additivity and prove Theorem \ref{main}. In Section 5 we give some discussions and problems.

\section{Norm defined on 1-homology group}
In article \cite{oo}, for an oriented, connected surface $K$ embedded in $S^{4}$, a norm is defined on $H_1(K,\mathbb{Z})$. Their construction applies to knotted tori. In this section, we introduce some first properties of a norm following their approach. 

We begin our discussion by defining a norm similar to the Thurston norm. For any oriented connected surface $F$, let the complexity of $F$, $x(F)$ be $max(-\chi (F), 0)$. If $F$ has more than one component, then $x(F)$ is the sum of them. 

For an arbitrary $n\geq 2$, let $K\colon K\to S^{n+2}$ be a codimension-two submanifold in $S^{n+2}$, namely a locally flat embedding from the manifold $K$ to the $(n+2)$-sphere. Denote the exterior of $K$ obtained by removing an open regular neighborhood of $K$ by $X_K$. 

To define the norm above, we need the lemma below:
\begin{lemma}
Let $K$ be a closed orientable $n$-manifold, and $Y$ be a simply connected closed (n+2)-manifold. Suppose $K\colon K\to Y$ is a null-homologous, locally flat embedding. Then $\partial X_K$ is canonically homeomorphic to $K \times S^1$, up to isotopy, such that the homomorphism $H_1(K)\to H_1(X_K)$ induced by including $K$ as the first factor $K \times pt$ is trivial, and for any slope $c\times pt$ there is a locally flat surface immersed in $X_K$ which bounds it.

\end{lemma}

The original lemma is Lemma 3.1 in article \cite{oo} for dimension 4, but in fact the proof is also true for high dimensions. We reformulate their proof as below:

\begin{proof}That $K$ is null-homologous induces that $K$ has a trivial normal bundle in $Y$, so $\partial X_K$ has a natural circle bundle structure $p\colon \partial X_K\to K$ over $K$. Framings of the normal bundle give the splitting of the bundle, and $H^1(X_K)\cong \mathbb{Z}$ and $H_1(X_K, \partial X_K)=0$ by Poincaré duality and excision. Thus, the homomorphism $H^1(X_K)\to H^1(\partial X_K)$ is injective, and the generator of $H^1(X_K)$ induces a homomorphism $\alpha\colon H_1(\partial X_K)\to \mathbb{Z}$. $\alpha$ sends the circle fiber of $\partial X_K$ to $\pm 1$, so the kernel of $\alpha$ projects isomorphically onto $H_1(K)$ via $p_\ast$. This shows $\partial X_K=K\times S^1$. It follows clearly from the construction that $H_1(K)\to H_1(X_K)$ is trivial. Moreover, if $c\times pt$ is an essential simple closed curve on $K\times pt$, it is homologically trivial in $X_K$, so it represents an element $[a_1, b_1]. . . [a_k, b_k]$ in $\pi_1(X_K)$. By a general position argument we may assume there is a surface to be a locally flat proper immersion that bounded by $c\times pt$.

\end{proof}

Therefore, for any homology class $r$ in $H_1(K,\mathbb{Z})$, we can define the complexity of $r$, $x(r)$ as the minimal possible $x(F)$ such that $F$ is a (possibly disconnected) oriented surface immersed in $X_K$ and its boundary expresses $r\times pt\in H_1(K\times S^1)$. Then $x(r_1+r_2)\leq x(r_1)+x(r_2)$ and $x(nr)\leq nx(r)$ because the possible bounded surfaces of $r_1+r_2$ include the union of $r_1$'s and $r_2$'s.

\begin{definition} $\|r\|_K=inf_n\frac{x(nr)}{n}$. 
\end{definition}

Now we prove it is actually a semi-norm. 

\begin{lemma}
(1)$\|nr\|_K=n\|r\|_K$. 

(2)$\|r_1+r_2\|_K\leq \|r_1\|_K+\|r_2\|_K$. 
\end{lemma}

\begin{proof}
For statement (1), because $x(nr)\leq nx(r)$, $\|nr\|_K\leq n\|r\|_K$. On the other hand, by definition $\|r\|_K=inf_m\frac{x(mr)}{m}\leq inf_m \frac{x(nmr)}{nm}=\frac{\|nr\|_K}{n}$. Therefore $\|nr\|_K=n\|r\|_K$. 

For statement (2), assume there are two surfaces $F_1$ and $F_2$ such that $\partial F_i=n_ir_i$ and $\|r_i\|_K\leq \frac{x(F_i)}{n_i}\leq \|r_i\|_K+\epsilon$, then $\partial (n_1F_2+n_2F_1)=n_1n_2(r_1+r_2)$ and $\frac{x(n_1F_2+n_2F_1)}{n_1n_2}=\frac{x(F_1)}{n_1}+\frac{x(F_2)}{n_2}\leq \|r_1\|_K+\|r_2\|_K+2\epsilon$, so $\|r_1+r_2\|_K\leq \|r_1\|_K+\|r_2\|_K+2\epsilon$, which means $\|r_1+r_2\|_K\leq\|r_1\|_1+\|r_2\|_2$. 

\end{proof}

This shows it is a semi-norm with $\mathbb{Z}$-coefficient. It can extend to $\mathbb{Q}$-coefficient by division, and $\mathbb{R}$-coefficient by continuity. In our proof $\mathbb{Z}$-coefficient is enough.

Generally speaking, we cannot judge whether $\|\cdot\|_K$ is degenerate, namely $\|x\|_K=0$ if and only if $x=0$ in $H_1(K,\mathbb{Z})$. We will prove that in some special case it is non-degenerate, and use it to prove theorem \ref{main}.

\begin{remark} To define a norm is a regular tool for proof. The norm we defined is similar to the most classical norm, Thurston norm, which is also first defined on $H_2(M^3)$ and then extend to $\mathbb{Q}$ and $\mathbb{R}$ coefficients: for $x\in H_2(M^3,\mathbb{Z})$, $\|x\|=min_Fx(F)$ where F runs over all properly embedded surfaces that represent $x$ (but don't have to be connected). An interesting problem of Thurston norm is how it can appear like, and a useful conclusion of Thurston norm is proved by Thurston: every symmetric integer polygon in $\mathbb{Z}^2$ with vertices satisfying the parity condition is the dual unit ball of the Thurston norm on a 3-manifold(\cite{thur}). 

Another useful norm is the intersection norm, it is defined with a set of loops for a surface M on $H_1(M,\mathbb{Z})$. It has some connection with the Thurston norm, see \cite{unit} for more details. 
\end{remark}

\section{Additivity of connected sum}

In this section we prove the key theorem. We use $\#$ to denote the connected sum. For an arbitrary $n\geq 2$, assume $Q=Q_1\#Q_2$ is a codimension-two submanifold embedded in $S^{n+2}$, where $Q_1$ and $Q_2$ are disjoint $n$-submanifolds embedded in $S^{n+2}$ with no link, namely they can be separated by an $(n+1)$-sphere standardly embedded in $S^{n+2}$.                                                                                                                                                                                                                                                                                                                                                                                                                                                                                                                                                                                                                                                                                                                                                                                                                                                                                                                                                                                                                                                                                                                                                                                                                                                                                                                                                                                                                                                                                                                                                                                                                                                                                                                                                                                                                                                                                                                                                                                                                                                                                                                                                                                                                                                                                                                                                                                                                                                                                                                                                                                                                                                                                                                                                                                                                                                                                                                                                                                                                                                        Then for any $r\in H_1(Q,\mathbb{Z})=H_1(Q_1,\mathbb{Z})\oplus H_1(Q_2,\mathbb{Z})$, $r=r|_{Q_1}+r|_{Q_2}$, we will prove the additivity of connected sum:

\begin{theorem}\label{key}
$\|r\|_Q=\|r|_{Q_1}\|_{Q_1}+\|r|_{Q_2}\|_{Q_2}$.
\end{theorem}

\begin{remark}In the knot theory, there is a similar equation $g(s_1\#s_2)=g(s_1)+g(s_2)$, where $s_1$ and $s_2$ are two knotted slopes embedded in $S^3$ and $g(s)$ is the genus (the minimal possible genus of surface bounded by $s$ embedded in $S^3$). We will also use the similar method to prove it.
\end{remark}

\begin{proof}
Denote the decomposition by $r=r_1+r_2$, and $\|\cdot\|_{Q_1}=\|\cdot\|_1$, $\|\cdot\|_{Q_2}=\|\cdot\|_2$. \\

$\|r\|_Q\leq \|r_1\|_1+\|r_2\|_2$ is obvious. Assume there are two surfaces $F_1$ and $F_2$ such that $\partial F_i=n_ir_i$ and $\|r_i\|_i\leq \frac{x(F_i)}{n_i}\leq \|r_i\|_i+\epsilon$, then $\partial (n_1F_2+n_2F_1)=n_1n_2r$ and $\frac{x(n_1F_2+n_2F_1)}{n_1n_2}\leq \frac{x(F_1)}{n_1}+\frac{x(F_2)}{n_2}\leq \|r_1\|_1+\|r_2\|_2+2\epsilon$, so $\|r\|_Q\leq \|r_1\|_1+\|r_2\|_2+2\epsilon$ which means $\|r\|_Q\leq\|r_1\|_1+\|r_2\|_2$. 

To prove $\|r\|_Q\geq \|r_1\|_1+\|r_2\|_2$, we need some surgeries.

\begin{figure}[h] 
		\centering
		\includegraphics[width=0.7\textwidth]{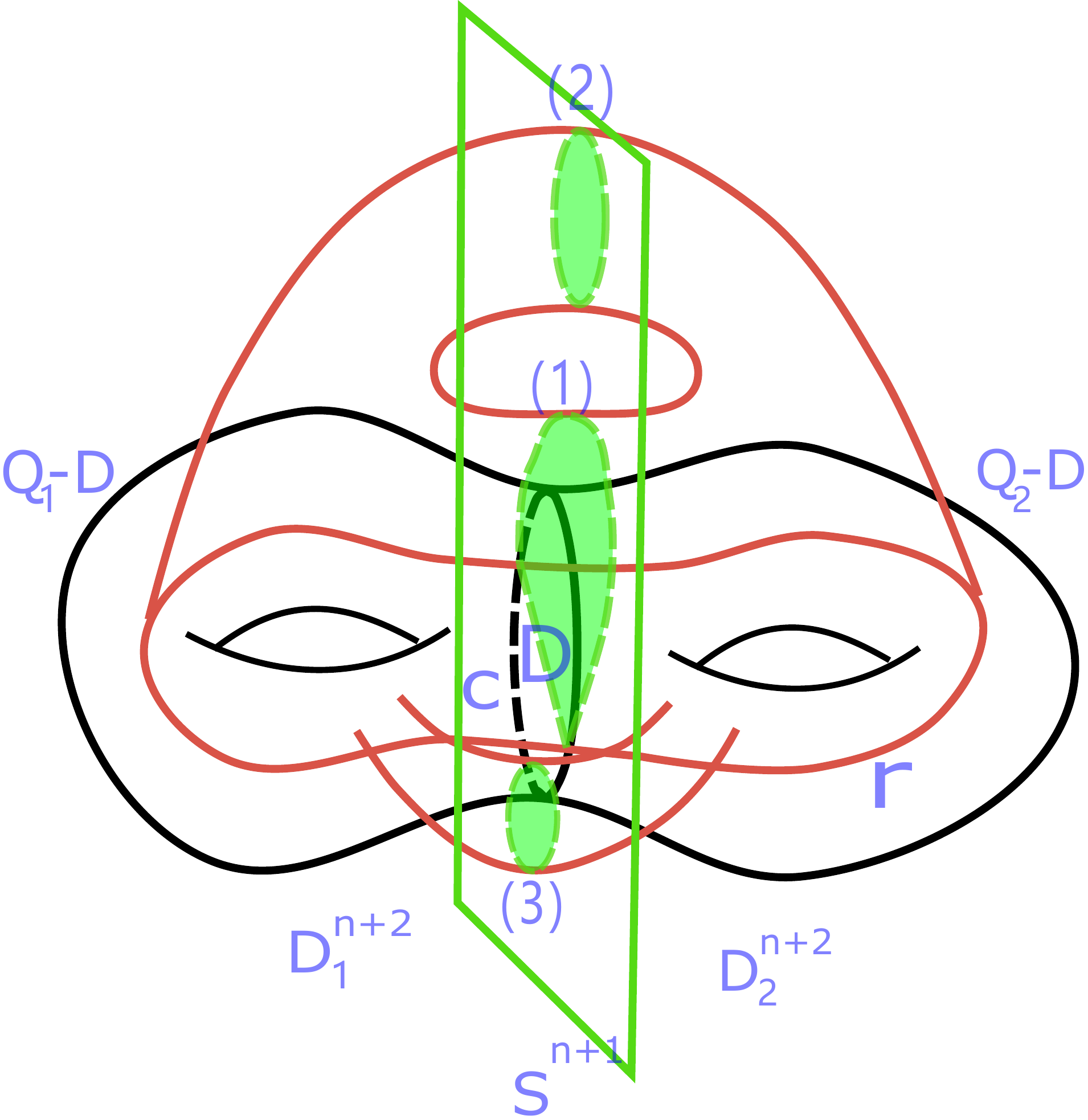} 
		\caption{surgery} 
		\label{Fig.main1}
	\end{figure}

Assume there is a surface $F$ smoothly immersed (disturb it if not) into $S^{n+2}$, such that $\partial F=nr$ and $\|r\|_Q\leq \frac{x(F)}{n}\leq \|r\|_Q+\epsilon$. As shown in the figure, we may use a $S^{n+1}$ to divide the connected sum and $S^{n+2}$ into $P_1=Q_1-D \subset D_1^{n+2}$ and $P_2=Q_2-D \subset D_2^{n+2}$, where $D_1^{n+2}\cup D_2^{n+2}=S^{n+2}$, $P_1\cup P_2=Q$, and $D$ is the connecting $n$-disk of $Q_1$ and $Q_2$ embedded in $S^{n+1}$. Denote the boundary of $D$ by $c=S^{n-1}$, which is trivially embedded in $S^{n+1}$. Then by Alexander duality, we see that $\pi_1(S^{n+1}-c)=H_1(S^{n+1}-c)=H_1(S^{n+2}-Q)=H_1(D_1^{n+2}-P_1)=H_1(D_2^{n+2}-P_2)=\mathbb{Z}$ has the same generator which is a circle $d$ linked with $c$ in $S^{n+1}$.

The intersection $F\cap S^{n+1}$ consists of curves in three cases:\\

{\it(1) curves with two endpoints on $c$. 

(2) closed curves not linked with $c$. (Namely they have 0 homology in $H_1(S^{n+1}-c)$).

(3) closed curves linked with $c$. (Namely they have non-zero homology in  $H_1(S^{n+1}-c)$).}\\

We first discuss the case (3).  Let the curves of case (3) have forms $m_1d, m_2d, . . . $ and $md$ be their common multiple. The generator in $H_1(S^{n+1}-c)=H_1(S^{n+2}-Q)=\mathbb{Z}$ is $d$, so the pre-image of $md$ under the map $\pi_1(F)\to \pi_1(S^{n+2}-Q)\to H_1(S^{n+2}-Q)$ has index $m$, which decides a covering of $F$, which we denote by $G$, and the covering degree is $m$. Besides, $G\to F\hookrightarrow (S^{n+2}-Q)$ gives an immersing from $G$ to the whole space, and the curves of $G$ cutting $S^{n+1}$ in case (3) represent the same homology element $w=md$ or its opposite in $H_1(S^{n+1}-c)=H_1(S^{n+2}-Q)=\mathbb{Z}$. We have $x(G)=mx(F)$ and $\partial G=m\partial F=mnr$, so $\frac{x(G)}{mn}=\frac{x(F)}{n}\leq \|r\|_Q+\epsilon$. Thus, we can use $G$ to replace $F$ and change (3) as\\

{\it (3)* closed curves linked with c, representing homology element $w$ or $-w$ in $H_1(S^{n+1}-c)=H_1(S^{n+2}-Q)=H_1(D_1^{n+2}-P_1)=H_1(D_2^{n+2}-P_2)=\mathbb{Z}$. }\\

For case (1), we cut $G$ along the curve and use an arc in $D$ to connect the endpoints and repair the surface by two $2$-disks immersed in $S^{n+1}$ and disjoint with $c$, bounded by the curve and the arc, so that we can change $\partial G$ to be in $Q_1$ and $Q_2$ representing $mnr_1$ and $mnr_2$. This surgery does not make $x(G)$ increase. 

For case (2), we cut $G$ along the curve and repair it by disks immersed in $S^{n+1}$ and disjoint with $c$. Also, $x(G)$ does not increase. 

For case (3)*, we note that $r_i\subset D_i^{n+2}-P_i$ represents 0 in $H_1(D_i^{n+2}-P_i)$, so the closed curves of $G$ cutting $S^{n+1}$ in the cases (2) and (3)* must have 0 homology class in $H_1(S^{n+1}-c)=H_1(D_i^{n+2}-P_i)$ in total. The case (2) is 0 homology, so case (3)* consists of $+w$ and $-w$ curve-pairs. For each pair, we cut $G$ along the two curves and try to connect the pair by a tube in $D_1^{n+2}-P_1$ and $D_2^{n+2}-P_2$. This is easy to obtain because we can connect $w$ and $-w$ in $H_1(S^{n+1}-c)=\pi_1(S^{n+1}-c)$ along an arc on the tubular neighborhood of $P_i$. By definition, this surgery does not change $x(G)$. \\

After doing the surgeries above, the original $G$ becomes two parts $G_1\subset D_1^{n+2}-(Q_1-D)$ and $G_2\subset D_2^{n+2}-(Q_2-D)$, $\partial G_i=mnr_i$, $x(G)\geq x(G_1+G_2)=x(G_1)+x(G_2)$. Thus\\

$\|r\|_Q+\epsilon\geq \frac{x(F)}{n}=\frac{x(G)}{mn}\geq \frac{x(G_1)+x(G_2)}{mn}=\frac{x(G_1)}{mn}+\frac{x(G_2)}{mn}\geq \|r_1\|_1+\|r_2\|_2$. \\

This proves $\|r\|_Q\geq \|r_1\|_1+\|r_2\|_2$. 
\end{proof}

\section{Construction of the Norm by Additivity}

In this section, we will apply the additivity for the case $n=2$. We first introduce a theorem of the norm in \cite[p 134]{oo}.

\begin{lemma}\label{base}
For some $a,b$ (which can respectively have infinitely many positive integer values), there is a $K=T^2$ embedded in $S^4$, such that $\|Ax+By\|_K=a|A|+b|B|$ where $x,y$ are the basis of $H_1(T^2)$. 
\end{lemma} 

Theorem \ref{key} in Section 3 will give a way to obtain a norm expression for high genus surface. By using the additivity of connected sum, and using lemma \ref{base}, we get

\begin{lemma}\label{combin}
For some $a_1,a_2,...a_{2g}$ (which can respectively have infinitely many positive integer values), there is a surface $K=F_g$ of genus g embedded in $S^4$, such that $\|A_1x_1+...+A_{2g}x_{2g}\|_K=a_1|A_1|+...a_{2g}|A_{2g}|$ where $x_1,...x_{2g}$ are the standard basis of $H_1(F_g)$.
\end{lemma}

Thus, we can prove a conclusion of high genus cases:\\
 
\noindent\textbf{Theorem 1.1.}\emph{
There is a surface $M=F_g$ of genus g embedded in $S^4$ such that, its image of extendable self-homeomorphisms in $\mathrm{Aut}(M,\mathbb{Z})=\mathrm{Sp}(2g,\mathbb{Z})$(with some symplectic basis) can only be diagonal matrixes with diagonal elements $\pm1$.}\\

\begin{proof}
Take different $A_i$ in lemma \ref{combin}. For a fixed extendable map, it keeps the norm, so its matrix can only be a diagonal matrix with diagonal elements $\pm1$.
\end{proof}

Especially, in the case $g=1$, $\mathrm{MCG}(T^2)=\mathrm{Aut}(T^2,\mathbb{Z})$ so there is an embedding from $T^2$ to $S^4$ such that its extendable self-homeomorphisms can only have matrixes $\pm I$.

\section{Further discussion}

The norm method defines a norm on $H_1(M,\mathbb{Z})$ so the best result we can obtain is in homological scale. We can only use it to detect the message about $E(i)$ in  $\mathrm{Aut}(F_g,\mathbb{Z}))=\mathrm{Sp}(2g,\mathbb{Z})$ but not in $\mathrm{MCG}(F_g)$ when $g>1$. Article \cite{oo} gives the result of $\mathrm{MCG}(T^2)$, which is because $\mathrm{MCG}(T^2)=\mathrm{Aut}(T^2,\mathbb{Z})$ when genus is 1. To obtain more information of $\mathrm{MCG}(F_g)$, we need more discussion about the Torelli group. Thus, a general conjecture is:

\begin{conjecture}There is some surface $F_g$ embedded in $S^4$ such that it has finitely many mapping classes of extendable self-homeomorphisms.
\end{conjecture}

Another possible problem is to consider the case of high dimension. 

\begin{conjecture}There is some $T^p$ embedded in $S^{p+2}$ such that it has finitely many mapping classes(in $\mathrm{MCG}(T^p)$) or linear mapping classes(in $\mathrm{SL}(p,\mathbb{Z})$) of extendable homeomorphisms.
\end{conjecture}

\section*{Acknowledgements}

I extend my sincere gratitude to my mentor, Liu Yi, and the institution that nurtured me, Peking University. Besides, I thank for the suggestions about this article and patience from the anonymous referee.

\end{document}